\documentclass[11pt]{article}
\usepackage{amssymb}
\newtheorem{thm}{Theorem}[section]
\newtheorem{defn}[thm]{Definition}

\newtheorem{conj}[thm]{Conjecture}
\newtheorem{constr}[thm]{Construction}

\newcounter{example}
\newenvironment{example}
  {\refstepcounter{example}\begin{trivlist}
    \item[\hskip\labelsep\textbf{Example \theexample}]}
  {\hfill$\blacksquare$\end{trivlist}}
\newenvironment{remark}
  {\refstepcounter{thm}\begin{trivlist}
    \item[\hskip\labelsep\textbf{Remark \thethm}]}
  {\hfill$\blacksquare$\end{trivlist}}
\newcommand{\ba}{\leftarrow}
\newcommand{\by}{\times}

\newcommand{\cl}{\widetilde}
\newcommand{\eq}{\sim}
\newcommand{\g}{{\mathfrak g}}
\newcommand{\ghat}{{\hat{\mathfrak g}}}

\newcommand{\inter}{\cap}

\renewcommand{\l}{\ell}
\newcommand{\M}{{\mathcal M}}

\newcommand{\pa}[1]{{\langle{#1}\rangle}}
\newcommand{\Plucker}{{Pl\"ucker}}
\newcommand{\pf}{\noindent \textbf{Proof.} }
\newcommand{\qed}{\hfill $\blacksquare$ \medskip}
\newcommand{\R}{{\mathcal R}}
\renewcommand{\sl}{{\mathfrak{sl}}}
\newcommand{\sumss}[2]{{\mathop{\sum_{#1}}\limits_{#2}}}
\newcommand{\SSYT}{{\mathop{\mathrm{SSYT}}\nolimits}}

\newcommand{\union}{\cup}
\newcommand{\Union}{\bigcup}
\newcommand{\w}{\omega}

\newcommand{\Z}{{\mathbb Z}}

\newcommand{\mybox}{{\begin{picture}(15,15)%
\put(0,0){\line(1,0){15}}\put(0,0){\line(0,-1){15}}%
\put(15,0){\line(0,-1){15}}\put(0,-15){\line(1,0){15}}\end{picture}}}
\newcommand{\row}[2]{\multiput(0,#2)(15,0){#1}{\mybox}}
\newcommand{\bp}{\begin{picture}} \newcommand{\ep}{\end{picture}}
\newcommand{\rb}[1]{\raisebox{#1\unitlength}}

\title{Pl\"ucker Relations on Schur Functions}
\author{Michael Kleber\thanks{Supported by an NSF Mathematical 
          Sciences Postdoctoral Research Fellowship.}\\
        Massachusetts Institute of Technology}
\date{}

\begin{document}
\maketitle

\begin{abstract}
We present a set of algebraic relations among Schur functions which are
a multi-time generalization of the ``discrete Hirota relations'' known
to hold among the Schur functions of rectangular partitions.  We prove
the relations as an application of a technique for turning \Plucker\ 
relations into statements about Schur functions and other objects
with similar definitions as determinants.  We also give a quantum 
analog of the relations which incorporates spectral parameters.
Our proofs are mostly algebraic, but the relations have a clear
combinatorial side, which we discuss.
\end{abstract}

\section{Introduction}
\label{sec_intro}

Consider the following relationship among the Schur functions
$s_\lambda$ where $\lambda$ is a rectangular partition:
\begin{equation}
\label{eq_basic}
s_\pa{m^\l} s_\pa{m^\l} = 
  s_\pa{m+1^\l} s_\pa{m-1^\l} + s_\pa{m^{\l+1}} s_\pa{m^{\l-1}}.
\end{equation}
Here $\pa{m^\l}$ is the partition with $\l$ parts each of size $m$,
whose Young diagram is an $\l\times m$ rectangle.
A.~N.~Kirillov noticed this fact as a relation among the
characters of finite-dimensional representations of $\sl_n$ while
studying the Bethe Ansatz for a one-dimensional system called the
generalized Heisenberg magnet~\cite{Ki}.

In later work, Kirillov and Reshetikhin observed that the relations
could be viewed as a discrete dynamical system, known to mathematical
physics as the discrete Hirota relations~\cite{KR}.  The initial
conditions are the characters of the fundamental representations of
$\sl_n$, and expressing the solutions in terms of the initial
conditions is precisely the Jacobi-Trudi formula for $s_\pa{m^\l}$.

In this paper, we present the natural extension of this set of
relations to Schur functions of arbitrary partitions.  The relations
are all of the form
\begin{equation}
\label{eq_otherterms}
s_\lambda s_\lambda = 
  s_{\lambda+\w_\l} s_{\lambda-\w_\l} + \mbox{ other terms}.
\end{equation}
Here we borrow notation from Lie theory: if $\lambda$ is a partition,
then we write $\lambda\pm\w_\l$ for the partition obtained by adding
or removing a column of height $\l$ from the Young diagram of
$\lambda$; this corresponds to taking the highest weight $\lambda$ and
adding or subtracting the fundamental weight $\w_\l$.  We have one
such relation for every choice of a partition $\lambda$ and column
height $\l$ such that $\lambda$ has a column of height $\l$ to begin
with (otherwise $\lambda-\w_\l$ does not make sense).  The various
choices of $\l$ should be thought of as independent time directions in
which we can evolve the dynamical system.

The ``other terms'' in equation~(\ref{eq_otherterms}) are also
products of two Schur functions, and all have coefficients $\pm1$.
The partitions that appear never have more columns or more outside
corners than $\lambda$ does.  Thus we get a hierarchy of sets of
relations for partitions with up to $k$ corners; when $k=1$ we are
restricted to rectangular partitions, and we recover
equation~(\ref{eq_basic}).

We prove the relations by reducing them to the \Plucker\ relations
among determinants of minors of a certain matrix, whose construction
we define in Section~\ref{sec_jt}.  The construction applies not only
to Schur functions, which we now view as determinants of the
Jacobi-Trudi matrices, but to the determinants of any family of
matrices with a similar type of definition.  We formalize this notion,
giving several other examples and a general version of the
construction.  The \Plucker\ relations themselves are briefly reviewed
in Section~\ref{sec_plucker}.

In Section~\ref{sec_main} we state and prove the relations.  We also
prove a generalization of the relations to ones which include
``shifts'' or ``spectral parameters.''  The generalizations of Schur
functions that satisfy this version of the equations are the quantum
analogs of characters for finite-dimensional representations for
$U_q(\widehat{\sl_n})$, and the generalized version may be related to
the representation theory of quantum affine algebras, which is not yet 
well understood.

Finally, while most of the earlier proofs are algebraic, in
Section~\ref{sec_combin} we offer a combinatorial interpretation for
the relations in terms of the Littlewood-Richardson rule, in which the
coefficients of $\pm1$ in the other terms mentioned above arise from
an inclusion-exclusion argument.  We give a completely bijective proof
for equation~(\ref{eq_basic}), the rectangular Young diagram version,
and we conjecture the existence of bijections with certain properties
that would lead to a fully combinatorial proof of
equation~(\ref{eq_otherterms}) as well.

The author is grateful to S.~Fomin and N.~Reshetikhin for helpful
discussions of the subject, and to W.~Brockman and S.~Billey for
comments on an earlier draft of this work.

\section{Generalized Jacobi-Trudi sets}
\label{sec_jt} \label{sec_def}

We will describe a scheme for translating the \Plucker\ relations
among determinants of minors of a matrix into identities of objects
defined by a Jacobi-Trudi style determinantal formula.  An instance of
this technique was used implicitly in \cite{LWZ} to prove some
relations among quantum transfer matrices.  Our applications will
include Schur functions (characters of representations of $SL_n$),
skew Schur functions, and Schur functions with spectral parameters
(quantum characters of $U_q(\widehat{\sl_n})$).  We construct the 
matrix in this section; the \Plucker\ relations are discussed in
Section~\ref{sec_plucker}.

The heart of the construction is an operation $(A,B) \to A\square B$,
where $A$ and $B$ are $n\by n$ matrices and $A\square B$ is an
$(n+1)\by(2n+2)$ matrix.  The operation can be depicted graphically
as:
\begin{equation}
\label{boxpicture}
\setlength{\unitlength}{1\baselineskip}
\raisebox{-1.5\unitlength}{\begin{picture}(3,3)
\put(0,0){\framebox(3,3){$A$}} \end{picture}}
\,\,,\,\,
\raisebox{-1.5\unitlength}{\begin{picture}(3,3)
\put(0,0){\framebox(3,3){$B$}} \end{picture}}
\,\,\to\,
\raisebox{-4\unitlength}{
\begin{picture}(4,8)
\put(0,0){\framebox(4,8){}}
\put(0,3){\line(1,0){4}}
\put(0,6){\line(1,0){4}}
\put(0,7){\line(1,0){4}}
\put(1,0){\line(0,1){3}}
\put(3,3){\line(0,1){3}}
\put(0,3){\makebox(3,3){$A$}}
\put(3,4){\makebox(1,1){$*$}}
\put(1,0){\makebox(3,3){$B$}}
\put(0,1){\makebox(1,1){$*$}}
\put(0,6){\makebox(1,1){$0$}}
\put(1,6){\makebox(1,1){$\cdots$}}
\put(2,6){\makebox(1,1){$0$}}
\put(3,6){\makebox(1,1){$\pm1$}}
\put(0,7){\makebox(1,1){$1$}}
\put(1,7){\makebox(1,1){$0$}}
\put(2,7){\makebox(1,1){$\cdots$}}
\put(3,7){\makebox(1,1){$0$}}
\end{picture}}
\end{equation}
We will first define the operation for our motivating example, the set 
of Jacobi-Trudi matrices:
$$
\left\{ M_\lambda := (h_{\lambda_i-i+j})_{i,j=1}^n \,|\,
 n\in\Z_{\geq0},\,\lambda\mbox{ a partition with $n$ parts}\right\} \!.
$$
If $h_k$ is the $k$th homogeneous symmetric function (so $h_0=1$ and
$h_k=0$ for $k<0$), then $\det(M_\lambda)$ is the Schur function
$s_\lambda$.  We permit $\lambda$ to end with zeros, so if $\lambda$
is a partition with $n$ parts then we can obtain $s_\lambda$ as the
determinant of such an $m\by m$ matrix for any $m\geq n$.

\begin{constr}
\label{constr_basic}
Let $\lambda$ and $\nu$ be partitions with $n$ parts.  We define the
matrix $M=M_\lambda\square M_\nu$, with $n+1$ columns indexed by
$\{1,\ldots,n+1\}$ and $2n+2$ rows indexed by
$\{L,R,1,\ldots,n,1',\ldots,n'\}$, as follows:
$$
\begin{array}{lll}
M_{Lj} &=& \delta(j,1) \\
M_{Rj} &=& (-1)^n\delta(j,n+1) \\
M_{ij} &=& h_{\lambda_i-i+j},\,\,i=1,\ldots,n \\
M_{i'j} &=& h_{\nu_i-i+j-1},\,\,i=1,\ldots,n
\end{array}
$$
\end{constr}
We adopt the notation $[r_1r_2\ldots r_k]_M$ for the determinant of
the $k\by k$ minor of a $k\by n$ matrix $M$ consisting of rows with
indices $r_1,\ldots,r_k$; when the choice of $M$ is clear from context
the subscript will be dropped.  Then for $M=M_\lambda\square M_\mu$,
we have $[R12\ldots n]=s_\lambda$ and $[L1'\ldots n']=s_\mu$.  (The
sign of $M_{R,n+1}$ was chosen for convenience precisely to make this
happen.)  \Plucker\ relations on $M$ will give us relations among
Schur functions.

The construction relies on the following property of the set of
Jacobi-Trudi matrices $\{M_\lambda\}$: there is a unique way to fill
in the $*$ regions in equation~(\ref{boxpicture}) so that any
$n+1$-row minor of $M_\lambda\square M_\mu$ with nonzero determinant
is again some $M_\nu$, possibly padded with one or two rows and
columns that do not affect the determinant.  To give a generalization
of the construction, we isolate this property.

\begin{defn}
\label{def_jt}
Let $\M$ be a set of square matrices.  Let $\R_n$ denote the set of
$n$-component vectors that appear as rows in any $n\by n$ matrix
$M\in\M$, for each $n\in\Z_+$.  We say $\M$ is a {\em generalized
Jacobi-Trudi set} if there exist equivalence relations $\eq_n$ on
$\R_n$ such that:
\begin{enumerate}
\item
Any two rows of an $n\by n$ matrix $M\in\M$ are $\eq_n$ related,
\item
If $M$ is an $n\by n$ matrix with nonzero determinant and all of its
rows are pairwise $\eq_n$ related, then there is a matrix $M'\in\M$
with the same rows as $M$ (but possibly permuted).
\end{enumerate}
Consider the operators $d_L$ and $d_R$, which respectively drop the
left and right components of a row vector.
\begin{enumerate} \addtocounter{enumi}{2}
\item
Take any two rows $r_1,r_2\in\R_n$ such that $d_L(r_1),
d_L(r_2)\in\R_{n-1}$.  If $r_1 \eq_n r_2$, then
$d_L(r_1)\eq_{n-1}d_L(r_2)$.  Furthermore, $d_L(r_1)=d_L(r_2)$ only if
$r_1=r_2$.  Thus we can talk about $d_L$ acting on the equivalence
classes.  Likewise, all this must hold for $d_R$ as well.
\item
If $A$ and $B$ are two $\eq_n$ classes such that $d_L(A)=d_R(B)$
then there is a unique $\eq_{n+1}$ class $C$ such that $d_R(C)=A$
and $d_L(C)=B$.
\end{enumerate}
\end{defn}

Our archetypical generalized Jacobi-Trudi set of matrices, of course,
is the set of Jacobi-Trudi matrices $M_\lambda$ defined above.  In
this case there is only one conjugacy class for each $\eq_n$, and it
consists of all rows of the form $(h_k,h_{k+1},\ldots,h_{k+n-1})$ for
$k+n-1$ nonnegative.  We should verify property (2) from the
definition: a matrix with $n$ rows of this form, ordered to have
nonincreasing values of $k$, corresponds to a partition unless some
row is repeated or some element on the diagonal is $h_j$ for $j$
negative.  In either case the resulting matrix has determinant 0.

Other examples of generalized Jacobi-Trudi sets include:

\begin{example}
\label{eg_formal}
The matrices $M_\lambda$ with the $h_k$ considered as formal
variables, {\em without} the specialization that $h_0=1$ and $h_k=0$
for $k<0$.  There is still only one $\eq_n$ class for each $n$, but
now every matrix formed of $n$ distinct rows from that class has
nonzero determinant.
\end{example}

\begin{example}
\label{eg_skew}
The matrices $M_{\lambda/\mu} := (h_{\lambda_i-\mu_j-i+j})_{i,j=1}^n$.
The determinants of these matrices are the skew Schur functions
$s_{\lambda/\mu}$ corresponding to skew Young diagrams $\lambda/\mu$,
with $\mu\subset\lambda$ ({\em i.e.}~$\mu_i\leq\lambda_i$ for all
$i$).  In this case, for each $n$ there are infinitely many $\eq_n$
classes, one for each choice of $\mu$: given a row vector
$(h_{a_1},h_{a_2},\ldots,h_{a_n})$, it can appear in matrices
$M_{\lambda/\mu}$ where $\mu_i-\mu_{i+1}=a_{i+1}-a_i-1$.

The operator $d_L$ (resp. $d_R$) takes the $\eq_n$ class
associated with $\mu$ to the $\eq_{n-1}$ class of $\mu$ with $\mu_1$
(resp. $\mu_{n-1}$) removed.  (Without loss of generality we
assume that $\mu_n=0$.)
\end{example}

\begin{example}
\label{eg_quantum}
The set of matrices $T_\lambda(u+c)$, where $\lambda$ is a partition,
$u$ is a formal variable, and $c\in\Z$ is called the shift.  We will
take the following as a formal definition:
$$
T_\lambda(u) := (t_{\lambda_i-i+j}(u+\lambda_1-\lambda_i+i+j-n-1))_{i,j=1}^n
$$
where $\lambda$ has $n$ parts, some of which may be zero.  Define
$s_\lambda^{(u)} := \det(T_\lambda(u))$.  As with Schur functions, the
$t_k(u)$ can optionally be specialized to $t_0(u)=1$, $t_k(u)=0$ for
$k<0$.

We will treat the $s_\lambda^{(u)}$ as formal symbols, but see the
remarks following Theorem~\ref{thm_quantum} for comments and
references on the mathematical physics origins of the objects.
Essentially, $s_\lambda^{(u)}$ can be regarded as quantum analogs of
characters of representations of $U_q(\ghat)$.  If we send the entry
$t_k(u+c)$ to $h_k$ and therefore ignore the shift (this is letting
$u\to\infty$ in the mathematical physics literature) we recover the
Jacobi-Trudi matrices $M_\lambda$ and plain Schur functions
$s_\lambda$.

To understand the equivalence classes here, note that the rows of any
matrix $T_\lambda(u+c)$ are of the form
$$(t_a(u+b),t_{a+1}(u+b+1),\ldots,t_{a+n-1}(u+b+n-1))$$
for some choice of integers $a$ and $b$.  The main diagonal of
$T_\lambda(u)$ has entries
$t_{\lambda_1}(*),t_{\lambda_2}(*),\ldots,t_{\lambda_n}(*)$,
while the anti-diagonal has
$t_*(u),t_*(u+\lambda_1-\lambda_2),\ldots,t_*(u+\lambda_1-\lambda_n)$.
It is therefore easy to see that if the row beginning with $t_a(u+b)$
appears in the matrix $T_\lambda(u+c)$, we must have $a+b =
\lambda_1-n+1$.  Therefore each $\eq_n$ class contains all rows which
share a common value $a+b$.

We remark that given a partition $\lambda$ with $n$ parts and an
$\eq_n$ class $A$, there is a unique integer $c$ such that the rows of
$T_\lambda(u+c)$ are in $A$.
\end{example}

Finally, we give a version of Construction~\ref{constr_basic} for any
generalized Jacobi-Trudi set of matrices, which we will apply to the
examples above.

\begin{constr}
\label{constr_general}
Let $\M$ be a generalized Jacobi-Trudi set of matrices, and take two
$n\by n$ matrices $A,B\in\M$.  Let $\cl{A},\cl{B}$ denote the $\eq_n$
classes of their respective rows.

We say $A$ and $B$ are {\em compatible} if
$d_L(\cl{A})=d_R(\cl{B})$.  For compatible $A,B$ we can
define the $(n+1)\by(2n+2)$ matrix $A\square B$.  Let $\cl{C}$ be the
$\eq_{n+1}$ class such that $d_R(\cl{C})=\cl{A}$ and
$d_L(\cl{C})=\cl{B}$, whose existence and uniqueness is
guaranteed by Definition~\ref{def_jt}. The rows of $A\square B$ are
indexed by $\{L,R,1,\ldots,n,1',\ldots,n'\}$.
\begin{itemize}
\item
Row $L$ is $(1,0,\ldots,0)$,
\item
Row $R$ is $(0,\ldots,0,(-1)^n)$,
\item
Row $i$ for $i=1,\ldots,n$ is the (unique) row $r_i\in\cl{C}$
such that $d_R(r_i)$ is the $i$th row of $A$,
\item
Row $i'$ for $i=1,\ldots,n$ is the (unique) row $r_{i'}\in\cl{C}$
such that $d_L(r_{i'})$ is the $i$th row of $B$.
\end{itemize}
\end{constr}
We point out again that $[R12\ldots n]=\det(A)$ and 
$[L1'\ldots n']=\det(B)$.  Of course, when $\M=\{M_\lambda\}$ this
reduces to Construction~\ref{constr_basic}.

\section{\Plucker\ relations}
\label{sec_plucker}

The \Plucker\ relations are the algebraic dependencies among the
determinants of the various minors of an arbitrary matrix.  We quickly
review them, fix our notation, and present a few pertinent examples.
For details on the subject, see \cite{Sturm}.

The \Plucker\ relations are most naturally stated for the $n\by n$
minors of an $n\by2n$ matrix.  Suppose we have such a matrix, and we
index its $2n$ rows by $1,\ldots,n,1',\ldots,n'$.  Pick some integer
$k$, $1\leq k\leq n$, and then pick $1\leq r_1<\cdots<r_k\leq n$.  The
relations state that
$$
[12\ldots n][1'2'\ldots n'] = 
  \sum_{1\leq s_1<\cdots<s_k\leq n}
  \sigma_{RS}([1,2,\ldots,n][1',2'\ldots,n'])
$$
where $\sigma_{RS}$ exchanges rows $r_i$ with $s_i'$ for
$i=1,\ldots,k$ before evaluating the determinants.  We could define
\Plucker\ relations more generally, for the minors of a matrix of any
size, but they would be a specialization of the above relations to the
case where some row in $[12\ldots n]$ is the same as a row in
$[1'2'\ldots n']$, resulting in many zeros on the right hand side of
the identity.

\begin{example}
\label{eg_7term}
The 7-term \Plucker\ relation arises from choosing $n=4$ and $k=2$.
It reads:
\begin{eqnarray*}
[{\bf1}{\bf2}34][5678] &=& 
  [5634][{\bf1}{\bf2}78]
+ [5734][{\bf1}6{\bf2}8]
+ [5834][{\bf1}67{\bf2}] +\\ & &
  [6734][5{\bf1}{\bf2}8]
+ [6834][5{\bf1}7{\bf2}]
+ [7834][56{\bf1}{\bf2}]
\end{eqnarray*}
In this example, we say rows $3$ and $4$ are {\em fixed}.
\end{example}

Naturally, we intend to apply the \Plucker\ relations to the
generalized Jacobi-Trudi sets of matrices defined in the previous
section.

\begin{example}
Consider the matrix
$$
M = \left[\begin{array}{cc}
h_3 & h_4 \\ h_2 & h_3 \\ h_1 & h_2 \\ 1 & h_1
\end{array}\right]
$$
Every $2\by2$ minor of $M$ is a Schur function.  Writing the Young
diagrams for their Schur functions, the 3-term \Plucker\ relation
$[{\bf1}2][34] = [42][3{\bf1}]+[32][{\bf1}4]$ yields
{
\setlength{\unitlength}{.75\unitlength}
$$
\rb{-15}{\bp(45,30) \row3{30} \row3{15} \ep} \,\,\,
\rb{-15}{\bp(15,30) \row1{30} \row1{15} \ep}
\, = \,
\rb{-15}{\bp(30,30) \row2{30} \row1{15} \ep} \,\,\,
\rb{-15}{\bp(45,30) \row3{30} \row2{15} \ep}
\, - \,
\rb{-15}{\bp(30,30) \row2{30} \row2{15} \ep} \,\,\,
\rb{-15}{\bp(45,30) \row3{30} \row1{15} \ep}
$$
}
The sign arises from the need to reorder the rows.
\end{example}

Finally, we will be looking at the \Plucker\ relations primarily for
the matrices $A\square B$ constructed in Section~\ref{sec_jt},
partitioning rows so that one term of the relation is
$\det(A)\det(B)$.  To specify an example of this type, we need to
choose matrices $A$ and $B$ from a generalized Jacobi-Trudi set, and
we need to pick some subset of the rows of either $A$ or $B$ (recall
that the $\square$ operation is not symmetric) to be the fixed rows in
the identity.

\begin{example}
\label{eg_321}
Take $\lambda=\pa{2,1,1}$ and $\mu=\pa{4,3,1}$.  We will apply the
7-term \Plucker\ relation (Example~\ref{eg_7term}) to
$T_\lambda(u)\square\, T_\mu(u)$ (Example~\ref{eg_quantum}).  Choosing
the first two rows of $T_\mu(u)$ as our fixed rows, we rearrange terms 
to get:
\begin{eqnarray*}
s_\pa{3,2,1}^{(u-1)} s_\pa{3,2,1}^{(u+1)}
&=&
 s_\pa{4,3,1}^{(u)} s_\pa{2,1,1}^{(u)} +
 s_\pa{3,2,2}^{(u-1)} s_\pa{3,1,1}^{(u+1)} +
 s_\pa{3,3,3}^{(u-1)} s_\pa{1,1,1}^{(u+3)}  \\ && {} +
 s_\pa{3,2,2,2}^{(u)} s_\pa{3,0}^{(u)} +
 s_\pa{3,3,3,2}^{(u)} s_\pa{1,0}^{(u+2)} -
 s_\pa{3,3,3,3}^{(u)} s_\pa{0,0}^{(u+3)}
\end{eqnarray*}
The order of terms was chosen as a precursor to Theorem~\ref{thm_main}.
If we ignore the spectral parameters, we get an identity on plain
Schur functions:

{
\setlength{\unitlength}{.5\unitlength}
\begin{eqnarray*}
\rb{-15}{\bp(45,45) \row3{45} \row2{30} \row1{15} \ep} \,\,\,
\rb{-15}{\bp(45,45) \row3{45} \row2{30} \row1{15} \ep}
&=&
\rb{-15}{\bp(60,45) \row4{45} \row3{30} \row1{15} \ep} \,\,\,
\rb{-15}{\bp(30,45) \row2{45} \row1{30} \row1{15} \ep}
\,+\,
\rb{-15}{\bp(45,45) \row3{45} \row2{30} \row2{15} \ep} \,\,\,
\rb{-15}{\bp(45,45) \row3{45} \row1{30} \row1{15} \ep}
\,+\,
\rb{-15}{\bp(45,45) \row3{45} \row3{30} \row3{15} \ep} \,\,\,
\rb{-15}{\bp(15,45) \row1{45} \row1{30} \row1{15} \ep}
\\[15\unitlength]
&& {} +\,
\rb{-30}{\bp(45,60) \row3{60} \row2{45} \row2{30} \row2{15} \ep} \,\,\,
\rb{15}{\bp(45,15) \row3{15} \ep}
\,+\,
\rb{-30}{\bp(45,60) \row3{60} \row3{45} \row3{30} \row2{15} \ep} \,\,\,
\rb{15}{\bp(15,15) \row1{15} \ep}
\,-\,
\rb{-30}{\bp(45,60) \row3{60} \row3{45} \row3{30} \row3{15} \ep}
\end{eqnarray*}
}

In the first version, the zero parts of the partitions are necessary
if the identity is to work without setting $t_0(u)=1$, $t_k(u)=0$ for
$k<0$.  If we are willing to make that specialization, we can drop the
zero parts, but we must adjust the shifts at the same time:
$s_\pa{\lambda_1\ldots\lambda_n,0}^{(u)} =
s_\pa{\lambda_1\ldots\lambda_n}^{(u-1)}$.  In the second version we
have already dropped the information about the zero parts.
\end{example}

\section{Main Theorem}
\label{sec_main}

In this section we present a set of recurrence relations, essentially
a discrete dynamical system, to which the Schur functions are a
solution.  These relations are a generalization of
equation~(\ref{eq_basic}), a system of relations which hold for the
Schur functions of partitions with rectangular Young diagrams.  We
also present the quantum analog of the relations, in
Theorem~\ref{thm_quantum} and following comments; this generalizes the
relation
\begin{equation}
\label{eq_qbasic}
s_\pa{m^\l}^{(u-1)} s_\pa{m^\l}^{(u+1)} = 
  s_\pa{m+1^\l}^{(u)} s_\pa{m-1^\l}^{(u)} + 
  s_\pa{m^{\l+1}}^{(u)} s_\pa{m^{\l-1}}^{(u)}
\end{equation}
We prove the relations by reducing them to \Plucker\ relations on
$M_\lambda\square M_\mu$, defined in Section~\ref{sec_jt}.  The simple
forms in equations~(\ref{eq_basic}) and~(\ref{eq_qbasic}) come from
the 3-term \Plucker\ relation.  Example~\ref{eg_321} in the previous
section shows the simplest instance based on the 7-term
\Plucker\ relation.

To state the relations, we first need to define some operations on the
partition $\lambda$, which we associate with its Young diagram
$Y=Y(\lambda)$.  Let $Y$ be a Young diagram with $n$ outside corners.
That is, we take $n$ points $(x_1,y_1),\ldots,(x_n,y_n)$ in
$\Z_{\geq0}\times\Z_{\geq0}$ with $x_1>\cdots>x_n$ and
$y_1<\cdots<y_n$, and the points in $Y$ are those less than any of the
$(x_i,y_i)$ in the product ordering.  We identify $Y$ with the
partition $\lambda=\pa{x_1^{y_1},x_2^{y_2-y_1},\ldots,x_n^{y_n-y_{n-1}}}$.
We also say that $Y$ has $n+1$ inside corners, numbered from $0$ to
$n$; the $i$th one has coordinates $(x_{i+1},y_i)$, where
$y_0=x_{n+1}=0$.

\begin{defn}
Let $Y$ be a Young diagram with $n$ outside corners as above, and pick 
two integers $i,j$ such that $1\leq i\leq j\leq n$.  We define two
Young diagrams by the coordinates of their corners:
$$
\begin{array}{l}
\pi^i_j(Y):
  \mbox{take the corners of $Y$, add $1$ to each of }
    x_{i+1},\ldots,x_j,
    y_{i},\ldots,y_j
\\
\mu^i_j(Y):
  \mbox{take the corners of $Y$, add $-1$ to each of }
    x_{i+1},\ldots,x_j,
    y_{i},\ldots,y_j
\end{array}
$$
These operations respectively add and remove a border strip which
reaches from the $i$th outside corner to the $j$th inside corner.

We will also want to add or remove several nested border strips.
Given integers $1\leq i_1<\cdots<i_r\leq j_r<\cdots<j_1\leq n$, we
further define
\begin{eqnarray*}
\pi^{i_1\cdots i_r}_{j_1\cdots j_r} &=& 
  \pi^{i_r}_{j_r}\circ\cdots\circ\pi^{i_1}_{j_1} \\
\mu^{i_1\cdots i_r}_{j_1\cdots j_r} &=& 
  \mu^{i_r}_{j_r}\circ\cdots\circ\mu^{i_1}_{j_1}
\end{eqnarray*}
Thus we add or remove border strips reaching from outside corner $i_s$
to inside corner $j_s$ for $1\leq s\leq r$.
\end{defn}

We apply these definitions of $\pi^{i_1\cdots i_r}_{j_1\cdots j_r}$
and $\mu^{i_1\cdots i_r}_{j_1\cdots j_r}$ only considering the
coordinates of corners, so the various $\pi^i_j$ and $\mu^i_j$
commute.  Note that applying $\pi^i_j$, for example, might decrease
the number of visible corners of $Y$ (by making $y_j$ the same as
$y_{j+1}$), but we ignore this effect in the latter definitions
above.  Since the intervals $[i_s,j_s]$ are nested, we will never end
up with $x_i<x_{i+1}$ or $y_i>y_{i+1}$.

Finally, we borrow notation from Lie theory: given a partition
$\lambda$, we let $\lambda\pm\w_\l$ denote the partition obtained from
$\lambda$ by adding or removing a column of height $\l$ to
$Y(\lambda)$.  If $\lambda=\pa{\lambda_1,\ldots,\lambda_m}$ and
$\mu=\lambda\pm\w_\l$, then $\mu_i=\lambda_i\pm1$ for $1\leq i\leq\l$
and $\mu_i=\lambda_i$ for $i>\l$.  Of course, we cannot take
$\lambda-\w_\l$ if $\lambda_\l=\lambda_{\l+1}$, that is, if
$Y(\lambda)$ does not have a column of height $\l$ to begin with.

\begin{thm}[Main Theorem]
\label{thm_main}
Take a partition $\lambda$ whose Young diagram $Y(\lambda)$ has $n$
outside corners.  Pick an integer $k$, $1\leq k\leq n$, and let $\l$
be the $k$th-shortest column height in $Y(\lambda)$, so $\l=y_k$ in
the coordinates above.  Then
$$
s_\lambda s_\lambda =
 s_{\lambda+\w_\l} s_{\lambda-\w_\l} +
 \sum_{r=1}^{\min(k,n-k+1)} \hspace{-1em}
  \sumss{1\leq i_1<\cdots<i_r\leq k}{k\leq j_r<\cdots<j_1\leq n}
   (-1)^{r-1}
   s_{\pi^{i_1\cdots i_r}_{j_1\cdots j_r}(\lambda)}
   s_{\mu^{i_1\cdots i_r}_{j_1\cdots j_r}(\lambda)}
$$
\end{thm}
That is, we take a signed double sum over all chains of properly
nested intervals $[i_1,j_1]\supset\cdots\supset[i_r,j_r]\ni k$.  For
each such chain we have the product of two Schur functions, obtained
by adding or removing all the corresponding border strips.

\begin{remark}
The recurrence relations can be viewed as defining the multi-time flow 
of a discrete dynamical system.  We think of $s_\lambda$ as being
associated with the lattice point whose $i$th coordinate is the number
of columns in $\lambda$ of height $i$.  If we allow arbitrary
partitions $\lambda$, the system is infinite-dimensional; if we
restrict ourselves to representations of $\sl_{n+1}$ it has dimension
$n$.

First, we note that that no partition appearing in
Theorem~\ref{thm_main} has more outside corners than $\lambda$ does.
Second, we observe that the only partition with more columns than
$\lambda$ is $\lambda+\w_\l$.  Therefore we can solve for
$s_{\lambda+\w_\l}$ to get a recurrence relation
$s_{\lambda+\w_\l}=(s_\lambda^2 - \sum\pm s_\pi s_\mu)/s_{\lambda-\w_\l}$,
expressing $s_{\lambda+\w_\l}$ in terms of Schur functions of
partition with strictly fewer columns and no more corners.  The
only initial conditions that need to be specified are for $s_\lambda$
when $\lambda$ has no two columns of the same height.
\end{remark}

\begin{example}
Take $\lambda$ to be the staircase partition $\pa{3,2,1}$ with $n=3$
corners, and pick $k=2$.  This instance of Theorem~\ref{thm_main} is
the Schur function part of Example~\ref{eg_321}.  The order in which
the terms appear there corresponds to taking the double sum over all
sets of nested intervals in the order:
$$
\underbrace{\{[2,2]\}\quad\{[1,2]\}\quad\{[2,3]\}\quad\{[1,3]\}}_{r=1}
\quad\underbrace{\{[1,3]\supset[2,2]\}}_{r=2}
$$
We will address the version with spectral parameters in
Theorem~\ref{thm_quantum}.
\end{example}

\pf
The formula is the \Plucker\ relation on $M_{\lambda-\w_\l} \square
M_{\lambda+\w_\l}$ in which we fix rows $1',\ldots,\l'$.  We index the
rows by $\{L,R,1,\ldots,m,1',\ldots,m'\}$ as in
Construction~\ref{constr_basic}, where $m$ is the number of parts of
$\lambda$.  The fixed rows are therefore those corresponding to rows
of $\lambda+\w_\l$ which got longer when the column of height $\l$ was
added.

First we locate the two pieces of Theorem~\ref{thm_main} outside the
double sum.  The term $s_{\lambda-\w_\l} s_{\lambda+\w_\l}$, of
course, is the \Plucker\ term $[R12\ldots m][L1'2'\ldots m']$, as we
have pointed out several times before.  The $s_\lambda s_\lambda$ term
is obtained from the \Plucker\ term $[L1\ldots\l(\l+1)'\ldots m']
[R1'\ldots\l'(\l+1)\ldots m]$, in which we swap $L$ with $R$ and every
row of $M_{\lambda+\w_\l}$ other than the fixed ones with the
corresponding row of $M_{\lambda-\w_\l}$.  This leaves rows $\l+1$
through $m$ of the two partitions unchanged in length.  The exchange
of $L$ and $R$ increases by one the lengths of rows $1$ through $\l$
of $\lambda-\w_\l$ and decreases by one the lengths of rows $1'$
through $\l'$ of $\lambda+\w_\l$, giving $\lambda$ in both cases.

To understand the remaining terms, we first observe that if the parts
of $\lambda$ are not all distinct, then the rows of $M_{\lambda-\w_\l}
\square M_{\lambda+\w_\l}$ are not all distinct either: if
$\lambda_i=\cdots=\lambda_j$ for some $i<j\leq\l$, then rows
$i,\ldots,j-1$ are exactly rows $(i+1)',\ldots,j'$.  Similarly, if
$\l<i<j$ then rows $i+1,\ldots,j$ are just rows $i',\ldots,(j-1)'$.
The nonzero terms in the \Plucker\ relation are determined by the
placement of the non-duplicated rows like $i'$ and $j$ (for $j\leq\l$)
or $i$ and $j'$ (for $i>\l$).  Thus we see that the number of corners
$n$, not the length $m$ of the partition, dictates the form of the
\Plucker\ relation we get.

Now consider the $[L1\ldots\l(\l+1)'\ldots m']
[R1'\ldots\l'(\l+1)\ldots m]$ term of the relation, already seen to
correspond to $s_\lambda s_\lambda$.  Observe that we can get any
other \Plucker\ term by exchanging some subset of $\{1,\ldots,\l\}$
from the first determinant with a subset of the same size drawn from
$\{R,(\l+1),\ldots,m\}$ from the second determinant.

Consider the effect of such a switch of a single row, say row $a$ for
row $b$, with $1\leq a\leq\l<b\leq m$.  As just noted, we may assume
$\lambda_a\neq\lambda_{a+1}$ and $\lambda_b\neq\lambda_{b-1}$.  Then
the effect on the $s_\lambda$ corresponding to
$[R1'\ldots\l'(\l+1)\ldots m]$ is as follows: we add an additional
part of size $\lambda_a$ to $\lambda$; we remove a part of size
$\lambda_b$; and we make each of $\lambda_{a+1},\ldots,\lambda_{b-1}$
one larger, corresponding to the newly shifted main diagonal of the
matrix $M_\lambda$.  But this is precisely the description of
$\pi^i_j(\lambda)$, where $Y(\lambda)$ has corner coordinates $y_i=a$
and $y_j=b-1$.  Likewise, in the $[L1\ldots\l(\l+1)'\ldots m']$ term
we change $\lambda$ into $\mu^i_j(\lambda)$.

One can easily check that the preceding argument still works when our
switched row $b$ is instead the one labeled $R$; in this case the
exchange has the effect of $\pi^i_n$ and $\mu^i_n$.  Exchanging
subsets larger than a single element is easily seen to mimic the
definition of $\pi^{i_1\cdots i_r}_{j_1\cdots j_r}$ and
$\mu^{i_1\cdots i_r}_{j_1\cdots j_r}$; the nesting of the intervals
arises because the ``push'' of $Y(\lambda)$ at outside corner $a$ and
the ``pull'' at inside corner $b$ are completely independent.

Finally, we find that after the swap of row $a$ for row $b$ described
above, we always need an odd number of adjacent transpositions to
correctly order the rows in the determinants: essentially, we need to
exchange rows $b$ and $b'$ but not $a$ and $a'$.  This explains the
$(-1)^{r-1}$ term, and completes the proof of Theorem~\ref{thm_main}.
\qed

There is a quantum analog of Theorem~\ref{thm_main} for the Schur
functions with spectral parameters defined in Example~\ref{eg_quantum}.

\begin{thm}
\label{thm_quantum}
For any partition $\lambda$, we can add spectral parameters to
the statement of Theorem~\ref{thm_main} to get
$$
s_\lambda^{(u-1)} s_\lambda^{(u+1)} =
 s_{\lambda+\w_\l}^{(u)} s_{\lambda-\w_\l}^{(u)} +
 \sum \sum \pm s_{\pi(\lambda)}^{(u+*)} s_{\mu(\lambda)}^{(u+*)}
$$
where the parameters inside the sum are as follows: given 
nested intervals $1\leq i_1<\cdots<j_1\leq n$, set
$\alpha=\pi^{i_1\cdots i_r}_{j_1\cdots j_r}(\lambda)$ and
$\beta=\mu^{i_1\cdots i_r}_{j_1\cdots j_r}(\lambda)$.  Then the
corresponding term in the sum is
$$
\begin{array}{rlll}
s_\alpha^{(u)}   & \hspace{-.5em} s_\beta^{(u+\lambda_1-\beta_1)} 
                     && \mbox{if } j_1=n \\[4pt]
s_\alpha^{(u-1)} & \hspace{-.5em} s_\beta^{(u+\lambda_1-\beta_1+1)} 
                     && \mbox{if } j_1<n
\end{array}
$$
\end{thm}
The case when $k=n$ is the subject of~\cite{LWZ}, where it is proved,
as here, by reducing to \Plucker\ relations.  Note that for $k=1$ or
$n$, the double sum is actually a single sum and no negative terms
appear.

\medskip\pf
As pointed out in Example~\ref{eg_quantum}, by appropriate choice of a
shift $c$, we can lift the matrix $M_\lambda$ to a matrix
$T_\lambda(u+c)$ whose rows are in whatever equivalence class we
choose.  Thus all we will do is pick some equivalence class, lift rows
$1,\ldots,m,1',\ldots,m'$ of $M_{\lambda-\w_\l} \square
M_{\lambda+\w_\l}$ to that class, and read off the necessary shifts
for each minor of our matrix to appear in the \Plucker\ relations.
Our choice of equivalence class is almost irrelevant; a different
choice would just correspond to adding a constant to $u$ in the final
relation.

We follow convention by choosing our equivalence class so that we are
dealing with minors of the matrix $M_{\lambda-\w_\l}(u) \square
M_{\lambda+\w_\l}(u)$, whose $2m$ rows other than $L$ and $R$ all look
like
$$
(t_{\lambda_1-c}(u-m+c), t_{\lambda_1+1-c}(u-m+1+c),\ldots,
 t_{\lambda_1+m-c}(u+c))
$$
The row with label $1'$ has this form with $c=0$, while the row with
label $1$ has $c=1$.  When we drop the left or right components of
these rows, respectively, we get the top rows of the matrices 
$M_{\lambda+\w_\l}(u)$ and $M_{\lambda-\w_\l}(u)$, as desired.

Given a minor corresponding to $s_*^{(u+c)}$, to identify the shift
$c$, recall that the top right entry in the matrix is $t_*(u+c)$.
Thus we can easily see that the $[R1'\ldots\l'(\l+1)\ldots m]
[L1\ldots\l(\l+1)'\ldots m']$ term of the \Plucker\ relation
corresponds to $s_\lambda^{(u-1)} s_\lambda^{(u+1)}$, again by looking
at the rows $1$ and $1'$ examined above.

Using the same reasoning, we see that for any $\alpha=\pi^{i_1\cdots
i_r}_{j_1\cdots j_r}(\lambda)$, the associated minor is either
$[R1'\ldots]$ (if row $R$ was not swapped away) or $[1'\ldots]$ (if
row $R$ was traded).  In the first case, we again end up with
$s_\alpha^{(u-1)}$; in the second case, we get $s_\alpha^{(u)}$.  Row
$R$ is swapped if and only if $j_1=n$, of course: this is the same as
saying the partition $\alpha$ has one more part than $\lambda$ if and
only if we added a border strip that reached the bottom row.

Determining the shift of $\beta=\mu^{i_1\cdots i_r}_{j_1\cdots
j_r}(\lambda)$ is more difficult because its top row, other than $L$
and possibly $R$, might be any of $1,2,\ldots,\l,\l+1$.  (Indeed, in
Example~\ref{eg_quantum}, each of these occurs.)  To sidestep this
difficulty, we note that the top row of the minor giving rise to
$\beta$ begins $t_{\beta_1}(*)$.  Assume that row $R$ was not traded.
Since we already know the top row must look like
$(t_{\lambda_1+1-c}(*),\ldots, t_{*}(u+c))$, we conclude that
$\beta_1=\lambda_1+1-c$, so $c=\lambda_1-\beta+1$.  Likewise, if row
$R$ was traded, the top row is one term shorter and ends with
$t_{*}(u+c-1)$, and the shift decreases by one, to
$\lambda_1-\beta_1$.
\qed

We conclude this section with a few comments on the relevance of the
quantum version of the theorem.

\begin{remark}
When we restrict $\lambda$ to being a partition with one corner, {\em
i.e.}~a rectangle, we are dealing with the 2-dimensional discrete
dynamical system
\begin{equation}
\label{eq_Q}
Q_{m+1}^\l(u) =
  \frac {Q_m^\l(u-1) Q_m^\l(u+1) - Q_m^{\l-1}(u) Q_m^{\l+1}(u)}
        {Q_{m-1}^\l(u)}
\end{equation}
for $\l=1,\ldots,n$ and $m\in\Z_+$.  Theorem~\ref{thm_quantum} states
that this system has a solution in which $Q_m^\l(u)$ is set to
$s_\pa{m^\l}^{(u)}$, an object which reduces to $s_\pa{m^\l}$ if we
ignore the spectral parameter.

The objects $s_\lambda^{(u)}$ themselves have a
representation-theoretic interpretation.  Using the fact that
$U_q(\widehat{\sl_n})$ is quasitriangular in the category of
evaluation modules and highest-weight modules, one can find a
commutative subalgebra which is a homomorphic image of the
Grothendieck ring.  The $s_\lambda^{(u)}$ live in this subalgebra and
play the role of characters; see~\cite{FR} for details.  These notions
come originally from mathematical physics and integrable systems,
where the objects in question are transfer matrices of Toda lattices;
see \cite{KNS},\cite{LWZ},\cite{BR} for more.

In this sense, dropping the spectral parameter corresponds to throwing
away some of the structure of the Lie algebra and retaining only the
action of the embedded subalgebra $U_q(\sl_n)$.
\end{remark}

\begin{remark}
Attempts to generalize this picture to Lie algebras of types other
than $A_n$ began in \cite{KR},\cite{KNS}.  In these cases, it appears
that the characters of $U_q(\ghat)$ do satisfy a generalized version
of equation~(\ref{eq_Q}).  Dropping the spectral parameters, though,
no longer gives statements about the fundamental representations of
$U_q(\g)$, but about certain non-irreducible representations which do
satisfy the discrete Hirota equations.  The generalization to Lie
algebras of types other than $A_n$ is hampered by the fact that there
is currently no character formula for representations of $U_q(\ghat)$.
\end{remark}

\begin{remark}
Recent work of the author~(\cite{Kl}) has shown a stronger result
about the generalized discrete Hirota relations, in an attempt to
sidestep the lack of a $U_q(\ghat)$ character formula.  For each Lie
algebra $\g$, there is a {\em unique} solution to the recurrence
relations in which $Q_m^\l$ is the character of a representation of
$U_q(\g)$ all of whose weights lie under $m\w_\l$ in the weight
lattice.  That is, we require that $Q_m^\l$ is a sum of irreducible
characters whose highest weights lie under $m\w_\l$, each occurring
with nonnegative integer coefficients.  This positivity constraint on
all of the infinitely many characters $Q_m^\l$ is quite rigid.

Theorem~\ref{thm_quantum} is the first step in extending this picture
from the rectangular case, in which only a small subset of the
characters of $U_q(\widehat{\sl_n})$ appear, to a full $n$-dimensional
system involving all of the characters.  One can hope that a
generalization of these new recurrence relations to other Lie algebras
may give us information on irreducible characters of $U_q(\ghat)$ for
which we do not even have conjectural values.
\end{remark}

\section{Combinatorial Considerations}
\label{sec_combin}

In this section, we look at the preceding formulas for Schur functions
purely combinatorially.  We offer a simple combinatorial proof of the
rectangle version of the formula, and indicate why we believe that the
subtraction that appears in Theorem~\ref{thm_main} arises from
inclusion-exclusion of sets labeled by single intervals.

We will multiply Schur functions using the following reformulation of
the Littlewood-Richardson rule, taken from \cite{Nak}, where the
technology of crystal bases is used to give an analog for Lie algebras
of type $B$, $C$, $D$ as well.

\begin{constr}
\label{constr_LR}
We wish to find the multiset $S$ of partitions such that $s_\lambda
s_\mu = \sum_{\nu\in S} s_\nu$.  To do this, let $\SSYT(\mu)$ be the
set of all semi-standard Young tableaux of shape $\mu$.  For any
tableau $T\in\SSYT(\mu)$, we obtain its column word
$cw(T)=i_1 i_2 \ldots i_m$ by reading off the numbers in $T$, reading each
column from top to bottom, beginning with the rightmost column and
ending with the leftmost.

Now we let the number $k$ act on the Young diagram $Y=Y(\lambda)$ by
adding one box to the $k$th row, provided $\lambda_k<\lambda_{k-1}$.
If $\lambda_k=\lambda_{k-1}$ then the action is illegal.  Denote the
resulting Young diagram by $Y\!\ba k$.  Then
$$
S = \left\{ (((Y\!\ba i_1)\ba i_2)\cdots\ba i_m) 
       \, \left|\strut\right. \,
    i_1 i_2 \ldots i_m = cw(T) \strut\right\}
$$
where $T$ ranges over all tableaux in $\SSYT(\mu)$ such that each
action is legal.
\end{constr}

Now we will give a purely bijective proof of the recurrence relation
for rectangular Young diagrams.  A proof was given in~\cite{Ki} which
did not mention the 3-term \Plucker\ relation, but which made use of
information from Lie theory about the dimensions of associated $\sl_n$
representations.

\begin{thm}
\quad \large  
$
s_\pa{m^\l} s_\pa{m^\l} = 
 s_\pa{m^{\l+1}} s_\pa{m^{\l-1}} + s_\pa{m+1^\l} s_\pa{m-1^\l}
$
\end{thm}

\pf 
Consider a tableau $T\in\SSYT(\pa{m^\l})$ such that the action of
$cw(T)$ on the Young diagram of shape $\pa{m^\l}$, as in
Construction~\ref{constr_LR}, is legal.  We consider two cases, based
on whether or not the leftmost column of $T$ consists exactly of the
numbers $1,2,\ldots,\l$.

If so, consider the tableau $T'$ obtained by removing the leftmost
column of $T$.  Observe that $T'\in\SSYT(\pa{m-1^\l})$, and the action
of $cw(T')$ on $Y(\pa{m+1^\l})$ is legal and yields the same Young
diagram as the action of $cw(T)$ on $Y(\pa{m^\l})$.  Furthermore,
all elements $T'$ of $\SSYT(\pa{m-1^\l})$ whose actions are legal
arise in this way; we need only note that $cw(T')$ never tries to
build on column $m+1$ of $Y(\pa{m+1^\l})$.

Otherwise, the leftmost column of $T$ contains an entry strictly
larger than $\l$, and therefore so does every column, as rows of $T$
are weakly increasing.  Now note that in any column of $T$, the
smallest number greater than $\l$ that appears must be $\l+1$.  This
is clear for the rightmost column, since $cw(T)$ acts legally on
$Y(\pa{m^\l})$, and can be seen inductively working to the left, again
because rows weakly increase.  Therefore we can consider the tableau
$T'$ obtained by removing the $\l+1$ from each column and pushing up
all the numbers below it; clearly $T'\in\SSYT(\pa{m^{\l-1}})$.  As in
the first case, this operation gives a bijection between $T$ acting
legally on $Y(\pa{m^\l})$ and $T'$ acting legally on
$Y(\pa{m^{\l+1}})$.
\qed

We are currently unable to provide a generalization of this argument
to arbitrary partitions $\lambda$, but we strongly believe that one
does exist.  Based on computational examples, we conjecture the
following form for a bijective proof of Theorem~\ref{thm_main}.

\begin{conj}
Let $\lambda$ be a partition with $n$ outside corners, choose a corner
$k$ and corresponding weight $\w_\l$, and retain the notions of
Theorem~\ref{thm_main}.  Let $L$ be the set of\/ $\SSYT(\lambda)$
acting legally on $Y(\lambda)$.
\begin{enumerate}
\item
The tableaux in $\SSYT(\lambda-\w_\l)$ which act legally on
$Y(\lambda+\w_\l)$ can be put in bijection with a subset $A$ of $L$.
\item
There are subsets $B^i_j\subseteq L\setminus A$, for each $1\leq i\leq
k\leq j\leq n$, such that $B^i_j$ is in bijection with
$\SSYT(\mu^i_j(\lambda))$ acting legally on $Y(\pi^i_j(\lambda))$.
\item
$L=A\union\Union B^i_j$.
\item
The intersection $B^{i_1}_{j_1}\inter\cdots\inter B^{i_r}_{j_r}$ is
nonempty if and only if we can reorder the terms to get
$1\leq i_1<\cdots<i_r\leq k\leq j_r<\cdots<j_1\leq n$, and in that
case it is in bijection with 
$\SSYT(\mu^{i_1\cdots i_r}_{j_1\cdots j_r}(\lambda))$ acting legally
on $Y(\pi^{i_1\cdots i_r}_{j_1\cdots j_r}(\lambda))$.
\end{enumerate}
All of the bijections between $\SSYT(\lambda)$ acting on
$Y(\lambda)$ and $\SSYT(\alpha)$ acting on $Y(\beta)$ should respect
the Young diagrams produced by the two actions.
\end{conj}

The conjecture implies Theorem~\ref{thm_main}, using
inclusion-exclusion to take the union $\Union B^i_j$.
We presently do not know the bijections or even how to identify the
sets $A$, $B^i_j$ in $L$.

\begin{example}
Taking the Schur function part of Example~\ref{eg_321} once again, the
only subtraction that takes place is of the term $s_\pa{3,3,3,3}
s_\pa{0,0}$, corresponding to the nested intervals
$[1,3]\supset[2,2]$.  There is one tableau (the empty tableau) whose
shape is the partition of zero.  To verify this instance of the
conjecture, we need to check that the Young diagram $Y(\pa{3,3,3,3})$
appears in the terms corresponding to intervals $[1,3]$ and $[2,2]$
once each.

This does happen: the element of $\SSYT(\pa{3,1,1})$ whose column word
is $44234$ acts on $Y(\pa{3,2,2})$, and the element of
$\SSYT(\pa{1,0})$ whose column word is $4$ acts on
$Y(\pa{3,3,3,2})$, both producing $Y(\pa{3,3,3,3})$.
\end{example}

\end{document}